\documentclass[12pt]{amsart}

\usepackage{amssymb}
\usepackage{amsmath}
\evensidemargin\oddsidemargin

\begin{document}

\setcounter{page}{1}

\newtheorem{theorem}{Theorem}
\newtheorem{lemma}[theorem]{Lemma}
\newtheorem{proposition}[theorem]{Proposition}
\newtheorem{remark}[theorem]{Remark}
\newtheorem{remarks}[theorem]{Remarks}
\newtheorem{THEO}{Theorem\!\!}
\newtheorem{DEFI}{Definition\!\!}
\renewcommand{\theTHEO}{}
\renewcommand{\theDEFI}{}

\newcommand{\eqnsection}{
\renewcommand{\theequation}{\thesection.\arabic{equation}}
    \makeatletter
    \csname  @addtoreset\endcsname{equation}{section}
    \makeatother}
\eqnsection

\def\AA{{\mathcal A}}
\def\BB{{\mathcal B}}
\def\c{{\mathbb C}}
\def\CC{{\mathcal C}}
\def\e{{\mathbb E}}
\def\ee{\mathrm{e}}
\def\EE{{\mathcal E}}
\def\eps{\varepsilon}
\def\FF{\mathcal{F}}
\def\HH{\mathcal{H}}
\def\ii{{\rm i}}
\def\II{\mathcal{I}}
\def\k{\kappa}
\def\LL{\mathcal{L}}
\def\lla{\lambda}
\def\n{{\mathbb N}}
\def\OO{{\mathcal O}}
\def\p{{\mathbb P}}
\def\pph{\varphi}
\def\r{{\mathbb R}}
\def\z{{\mathbb Z}}
\def\MM{\mathcal{M}}
\def\NN{\mathcal{N}}
\def\SS{\mathcal{S}}
\def\VV{\mathcal{V}}
\def\WW{\mathcal{W}}
\def\Xde{{\tilde X}}
\def\Ydee{{\tilde Y}^\eta}

\def\lacc{\left\{}
\def\lcr{\left[}
\def\lpa{\left(}
\def\lva{\left|}
\def\racc{\right\}}
\def\rcr{\right]}
\def\rpa{\right)}
\def\rva{\right|}

\def\Un{{\bf 1}}
\def\d{\, \mathrm{d}}
\def\qed{\hfill$\square$}
\def\elaw{\stackrel{d}{=}}

\title[Absolute continuity of Ornstein-Uhlenbeck process]
      {On the absolute continuity of multidimensional Ornstein-Uhlenbeck processes}

\author[Thomas Simon]{Thomas Simon}

\address{Laboratoire Paul Painlev\'e, U. F. R. de Math\'ematiques, Universit\'e de Lille 1, F-59655 Villeneuve d'Ascq Cedex. {\em Email} : {\tt simon@math.univ-lille1.fr}}



\begin{abstract} Let $X$ be a $n$-dimensional Ornstein-Uhlenbeck process, solution of the S.D.E.
$$\d X_t\; =\; AX_t  \d t \; +\; \d B_t$$
where $A$ is a real $n\times n$ matrix and $B$ a L\'evy process without Gaussian part. We show that when $A$ is non-singular, the law of $X_1$ is absolutely continuous in $\r^n$ if and only if the jumping measure of $B$ fulfils a certain geometric condition with respect to $A,$ which we call the exhaustion property. This optimal criterion is much weaker than for the background driving L\'evy process $B$, which might be very singular and sometimes even have a one-dimensional discrete jumping measure. It also solves a difficult problem for a certain class of multivariate Non-Gaussian infinitely divisible distributions. 
\end{abstract}

\maketitle

\section{Introduction and statement of the result}

If $Z$ is a real L\'evy process without Gaussian part, finding a necessary and sufficient condition on its jumping measure $\nu$ for the absolute continuity of $Z_t$ at some given $t>0,$ is a hard problem for which no sensible conjecture has been formulated as yet. One of the main difficulties for this formulation stems from the time-dependency of the absolute continuity property: if $\nu$ is infinite and has discrete support, then there are some situations where e.g. $Z_1$ is singular and $Z_2$ is absolutely continous. We refer to \cite{Wa} and Chapter 27 in \cite{S} for more on this topic as well as further references. When $Z$ is multidimensional, the problem becomes increasingly complicated and some partial results had been given in \cite{Ya}, involving conditions of geometrical nature on $\nu$.

On the other hand, the problem of absolute continuity may well become simpler, and yield weaker conditions, when considering certain functionals of $Z$. In the real case for example, one can show that

\begin{equation}
\label{lifsh}
\int_0^1\!Z_t \d t\; {\rm is \; a.c.}\;\; \Longleftrightarrow\;\; \nu \;{\rm is\; infinite.}
\end {equation}
Notice that the condition on the right-hand side is only equivalent to the non-atomicity of $Z_1$ - see Theorem
27.4 in \cite{S}. With a view towards the methods developed later in the present paper, let us give a short proof of the reverse inclusion in (\ref{lifsh}), the direct one being straightforward. If $\nu$ is infinite and $T^\eta_1, T^\eta_2$ denote the two first jumping times of $Z$ into $[-\eta, \eta]^c,$ then $\p[T^\eta_2 \ge 1]\to 0$ as $\eta\to 0.$ Besides on $\{T^\eta_2 < 1\}$ one can write 
$$\int_0^1\! Z_t\, \d t\; = \; (T^\eta_2 - T^\eta_1)\Delta Z_{T^\eta_1} \; +\; (1 - T^\eta_2)\Delta Z_{T^\eta_1} \; +\;\int_0^1 \lpa Z_t - \Un_{\{T^\eta_1 \le t\}}\Delta Z_{T^\eta_1}\rpa \d t$$
and the right-hand side is a.c. on $\{T^\eta_2 < 1\}$ for every $\eta >0,$ since conditionnally on $\FF^\eta$ which is the $\sigma$-field generated by all the information given by $Z$ on $[0,1]$ except $T^\eta_1,$ the variable $(T^\eta_2 - T^\eta_1)$ has uniform hence absolutely continuous law on $[0,T^\eta_2],$ the variable $\Delta Z_{T^\eta_1}$ is non-zero a.s. and $\FF^\eta$-measurable, and the remaining terms are $\FF^\eta$-measurable. 

This conditioning method together with, roughly speaking, a derivation procedure along certain jumping times, had been systematically developed in the monograph \cite{DLS} where various absolute continuity results and smoother properties were established for several functionals of Wiener and other L\'evy processes, such as $L_p$-norms, proper integrals along a given function, one-sided and two-sided suprema. In \cite{NS, K}, it was also applied to a class of real stochastic equations with non-linear drift driven by $Z$.

In this paper, we will deal with Non-Gaussian multidimensional Ornstein-Uhlenbeck processes, which are solutions to the S.D.E.
\begin{equation}
\label{OU2}
\d X_t\; =\; AX_t  \d t \; +\; B \d Z_t
\end{equation}
where $A$ is a real $n\times n$ matrix, $B$ a real $n\times d$ matrix and $Z$ a $d-$dimensional L\'evy process without Gaussian part. Adaptating without difficulty the discussion made in \cite{S} pp. 104-105 to the case where $A$ is not necessarily a multiple of the identity matrix, we see that the solution to (\ref{OU2}) is given in terms of some L\'evy integral:
\begin{equation}
\label{levy}
X_t \; = \; \ee^{tA}x\; +\; \int_0^t \ee^{(t-s)A}B\d Z_s, \quad t\ge 0,
\end{equation}
where $x\in\r^n$ is the initial condition. Apart from a natural extension of the Langevin equation, these processes have their roots in distribution theory because in the ergodic case their limit laws are known \cite{SY} to be operator-self decomposable   - see \cite{SWYY} for further results. Nowadays, Non-Gaussian OU processes are also quite popular in modelling \cite{BS}.

To state our result, we need some notation. We will assume that the reader is familiar with basic properties of L\'evy processes and jumping measures which can be found at the beginning of the two monographs \cite{B, S}. Setting $\{B_t, \; t\ge 0\}$ for the $\r^n$-valued L\'evy process $\{BZ_t, \; t\ge 0\}$ and $\nu^B$ for its L\'evy measure, we introduce the vector spaces
$$\BB_t\; =\; {\rm Vect} [\Delta B_s, \; s\le t]\quad {\rm and} \quad \AA_t\; =\;<A,\BB_t>, \quad t > 0,$$
where here and throughout, for every vector subspace $\EE\subset\r^n$ with basis $\{e_1,\ldots, e_p\},$ we use the notation
$$<A,\EE>\; =\; {\rm Vect} [A^{i-1} e_j, \; i = 1\ldots n, j = 1\ldots p].$$
Notice that actually $\AA_t = {\rm Vect} [A^{i-1} \Delta B_s, \; i = 1\ldots q, s\le t],$ where $q$ stands for the degree of the minimal polynomial of $A.$ Setting $\BB = {\rm Im} \,B\subset\r^n,$ the condition $<A,\BB>\, =\r^n$ or, in an equivalent matrix formulation,
\begin{equation}
\label{Rank}
{\rm Rank} \lcr B, AB,\ldots, A^{q-1}B\rcr\; =\; n,
\end{equation}
is well-known as a controllability condition on the deterministic linear system
$$x'_t\; =\; Ax_t\; +\; Bu_t$$
where $\{u_t, \; t\ge 0\}$ is some input function - see e.g. Chapter 1 in \cite{W} and the references therein. When (\ref{Rank}) holds, we will say that $(A,B)$ is controllable.

Let $\k$ denote the cyclic index of $A$, which is the maximal dimension of its proper subspaces. The condition $\k =1$ means that $A$ is a cyclic matrix i.e. there exists a generating vector $b\in\r^n$ such that $(b, Ab, \ldots, A^{n-1}b)$ forms a basis of $\r^n.$ This is also  equivalent to $q,$ that is the minimal polynomial of $A$ is in fact its characteristic polynomial. When $\k >1,$ it is the number of the invariant factors of $A$, viz. the unique number of subspaces $\AA_i\subset \r^n$ such that $\AA_1\oplus\ldots\oplus \AA_\k =\r^n,$ with each $\AA_i$ stable by $A$ and each $A/\AA_i$ a cyclic operator whose minimal polynomial $\alpha_i$ divides $\alpha_{i-1},$ $\alpha_1$ being the minimal polynomial of $A$ itself - see Chapter 4 in \cite{G} or Section 0.10 in \cite{W} for more precisions. 

Let $m = {\rm Dim}\, \BB = {\rm Rank}\, B.$ When (\ref{Rank}) holds, a result of M.~Heymann - see Theorem 1.2. in \cite{W} - entails that necessarily $m \ge\kappa.$ More precisely, there exist $\k$ linearly independent vectors $b_1, \ldots, b_\k\in\BB$ such that
\begin{equation}
\label{won}
\BB^1\; +\; \ldots\; +\; \BB^\k\; =\; \r^n
\end{equation}
with the notation $\BB^i = {\rm Vect} \lcr b_i, Ab_i,\ldots, A^{q-1}b_i\rcr$ for $i = 1\ldots \k.$ Actually, Heymann's result was originally more precisely stated, connecting the subspaces $\BB^i$ to the above $\kappa$ invariant factors of $A.$ Nevertheless we shall not need this in the sequel. Assuming (\ref{Rank}), a sequence $(b_1, \ldots, b_r)\in\BB$ of linearly independent vectors such that $\BB^1\, +\, \ldots\, +\, \BB^r\, =\, \r^n$ with the above notations will be called a {\em generating sequence} of $\r^n$ with respect to $(A,B),$ or simply a generating sequence when no confusion is possible. Notice that $r\le m.$ Besides, the very definition of $\k$ entails that necessarily $r \ge \k$ as well  - see the proof of Theorem 1.2 in \cite{W} for details. We now come to the central definition of this work:

\begin{DEFI} With the above notations, the L\'evy measure $\nu^B$ is said to
  {\em exhaust} $\r^n$ {\em with respect to} $A$ if $<A,\BB>\; =\r^n$ and there exists $r\in [\k, m]$ and a subspace $\HH_r\subset\BB$ of dimension $r$ such that $<A,\HH_r>\; =\r^n$ and $\nu^B(\HH_r\cap\HH^c) =+\infty$ for every hyperplane $\HH\subset\HH_r.$
\end{DEFI}

This definition is related to the conditions given in \cite{Ya} for the absolute continuity of multivariate infinitely divisible distributions, but it is less stringent since no arcwise absolute continuity is required, and since $\nu^B$ may be carried by any subspace with dimension $r\in [\k, m]$. Here, however, the important fact is that this subspace must be chosen with respect to $A$. Introducing finally the stopping time
$$\tau\; =\; \inf\{ t>0, \; \AA_t \, =\r^n\},$$
our result reads as follows:
\begin{THEO} If $A$ is non-singular, then one has
$$X_1\; {\rm is \; a.c.}\;\; \Longleftrightarrow\;\; \tau\, =\, 0\;\, {\rm a.s.}\;\; \Longleftrightarrow\;\;  \nu^B \;{\rm exhausts}\; \r^n\; {\rm w.r.t.}\; A.$$
\end{THEO}
Notice that the above equivalences are time-independent, so that when $X_1$ is a.c. then $X_t$ is a.c. as well for every $t >0.$ In other words absolute continuity is not a temporal property for Non-Gaussian OU processes, unlike L\'evy processes without Gaussian part. Besides, when $X_1$ is not a.c. then the first equivalence entails that $X_1$ is valued with positive probability in some fixed affine hyperplane of $\r^n,$ or equivalently that a certain one-dimensional projection of $X_1$ must have an atom. This, again, contrasts with L\'evy processes since $Z_1$ may be non-atomic and not absolutely continuous - see Theorem 27.19 in \cite{S}.

We stress that the variable $X_1$ itself is infinitely divisible, i.e. it is the distribution at time 1 of some L\'evy process $\{Y_t, \; t\ge 0\}$ valued in $\r^n$. Indeed, a straightfoward extension of Lemma 17.1 in \cite{S} shows that $X_1$ is ID without Gaussian part and L\'evy measure
$$\nu^X(\Lambda)\; =\; \int_{\r^n}\nu^B(dx)\int_0^1\Un_\Lambda (\ee^{sA} x) \d s, \quad \Lambda\in\BB(\r^n).$$
Hence, our result yields an optimal criterion of absolute continuity for a certain subclass of multivariate Non-Gaussian ID distributions, which we may call the {\em OU class}. To this end, one can check that if $\nu^X$ satisfies any condition given in \cite{Ya}, then $\nu^B$ exhausts $\r^n$ w.r.t. $A$, but that the converse inclusion is not necessarily true since $\nu^X$ may not be arcwise absolutely continous when $\nu^B$ exhausts $\r^n$ w.r.t. $A.$

As we mentioned before, the variables $X_t$ converge in law when $t\to\infty$ to some operator self-decomposable or OL distribution in $\r^n$ under an ergodicity assumption, that is when the eigenvalues of $A$ have all negative real parts and $\nu^B$ is log-integrable at infinity \cite{SY}. If in addition $\nu^B$ exhausts $\r^n$ w.r.t. $A,$ then it is easy to see that the limit distribution is also genuinely $n$-dimensional. Let us notice that the absolute continuity of non-degenerated OL distributions had been established in \cite{Y}. This is probably related to our result, even though absolute continuity and non absolute continuity properties are barely  stable under weak convergence. To make a true connection, one would need a stronger type of convergence such as convergence in total variation \cite{DLS}, but no such result seems available in the literature.
 
It follows from our definition that when $\nu^B$ exhausts $\r^n,$ then $(A,B)$ is necessarily controllable. On the other hand, when $(A,B)$ is not controllable then $\tau = +\infty$ a.s. so that from our result $X_1$ is not a.c. In a recent paper \cite{PZ} which was the starting point of this work, it was  proved that $X_1$ is absolutely continuous as soon as $(A,B)$ is controllable and the jumping measure $\nu$ of $Z$ is absolutely continuous in an open neighbourhood of 0. These conditions entail of course that $\nu^B$ exhausts $\r^n,$ but they are highly non-equivalent: for example when $A$ is cyclic and non-singular, then our result entails that $Z$ might be genuinely one-dimensional with an infinite, and possibly discrete,  jumping measure carried by some line in $B^{-1}(b)$ where $b$ is a generating vector of $A$, nevertheless $X_1$ will be a.c. in $\r^n.$ Our method is also very different from \cite{PZ}, which hinges upon a certain derivation procedure, made possible by the absolute continuity condition on $\nu$, along the jumping {\em sizes} of $Z$. Here, as shortly suggested above, we will differentiate along suitably chosen jumping times, the price to pay being the non-singularity assumption on $A$. Our time-derivation procedure is close to the one developed in \cite{K}, whose Theorem 1.1. actually entails that $X_1$ is a.c. when $A$ is non-singular and $\nu^B(\HH^c) = \infty$ for every hyperplane $\HH\subset\r^n$. But again this latter assumption is more stringent than our exhaustion property, as well as, to the best of our knowledge, all conditions given in the literature on Malliavin's calculus for jump processes - see the references given in \cite{NS, K, PZ}.
 
In the third section of this paper we will briefly describe what happens when the driving process $Z$ has some Gaussian component. By independence and by the linearity of the equation (\ref{OU2}), we get an analogous result which
only requires a small modification of the proof in the Non-Gaussian case. Then we will discuss a few examples in order to provide more geometrical insight on the exhaustion property. In particular we will give a complete description in the case $n=2.$ Finally, we will mention some counterexamples and open questions. Before this we now proceed to the 

\section{Proof of the theorem}

\subsection{Proof that $X_1$ is a.c. $\Rightarrow\;\tau = 0$ a.s.} Setting 
$$\BB_0\; =\; \bigcap_{t>0}\BB_t$$
which is a deterministic subspace of $\BB$ by the 0-1 law, it follows from the definition of $\BB_t$ and that of a jumping measure that necessarily $\nu(\BB_0^c) < +\infty.$ In particular, $\p[T > 1] > 0$ where $T = \inf\{ t>0, \; \Delta B_t \in\BB_0^c\}.$ Endow $\BB$ with a canonical Euclidean structure and set $\BB_0^\perp$ for the orthogonal supplementary of $\BB_0$ in $\BB$, which may be reduced to $\{0\}$ if $\BB_0=\BB.$ Decomposing the L\'evy process $\{B_t, \; t >0\}$ along the orthogonal sum $\BB = \BB_0\oplus\BB_0^\perp:$ 
$$B_t\; =\; B^0_t\; +\; B^\perp_t, \quad t >0,$$
notice that the L\'evy process $\{B^\perp_t, \; t >0\}$ is either the zero process or a compound Poisson process with drift coefficient, say, $b^\perp$. Hence, on $\{T > 1\},$ we deduce from (\ref{levy}) that 
$$X_1 \; = \; \ee^Ax\; +\; \int_0^1 \ee^{(1-s)A}\d B^0_s\; -\; \int_0^1 \ee^{(1-s)A}b^\perp \d s$$
where we use the notation $b^\perp = 0$ if $B^\perp\equiv 0.$ Last, writing for every $t\in\r$ 
\begin{equation}
\label{minime}
\ee^{tA}\; =\; \sum_{k=1}^q \psi_r(t) A^{r-1},
\end{equation}
where the $\psi_r$'s are certain real functions whose exact expression is given e.g. in \cite{G} Chapter 5, we see that
$$\int_0^1 \ee^{(1-s)A}\d B^0_s\; \in\; <A, \BB_0>\quad {\rm a.s.}$$
Putting everything together entails that if $X_1$ is a.c. then necessarily $<A, \BB_0>\, =\r^n.$ But from the definitions of $\BB_0$ and $\tau,$ this yields $\tau = 0$ a.s. 

\qed

\begin{remark} {\em The above proof shows also that if $\p[\tau > 0] > 0,$ then $X_1$ is valued in some fixed affine hyperplane with probability $\p[T>1].$ Notice that if $\BB_0 =\BB,$ then this probability is 1 and $(A,B)$ is not controllable, which entails that actually $\tau =\infty$ a.s.} 
\end{remark}

\subsection{Proof that $\tau = 0$ a.s. $\Rightarrow\, \nu^B$ exhausts $\r^n$ w.r.t. $A$} We may assume $<A,\BB>\, =\r^n$ since otherwise $\tau =\infty$ a.s. Suppose now that $\nu^B$ does not exhaust $\r^n.$ By the above assumption there exists a hyperplane $\HH_1\subset\BB$ such that $\nu^B(\HH_1^c) < \infty.$ If $<A,\HH_1>\, \neq\r^n,$ then
$$\tau\; \ge \; \inf\{t > 0, \; \Delta B_t \in \HH_1^c\}\; >\; 0 \quad{\rm a.s.}$$
If $<A,\HH_1>\, =\r^n,$ then there exists a hyperplane $\HH_2\subset\HH_1$ such that $\nu^B(\HH_1\cap\HH_2^c) < \infty$ because $\nu^B$ does not exhaust $\r^n,$ and in particular $\nu^B(\HH_2^c) < \infty$ since we already have $\nu^B(\HH_1^c) < \infty.$ Similarly, when $<A,\HH_2>\, \neq\r^n,$ then
$$\tau\; \ge \; \inf\{t > 0, \; \Delta B_t \in \HH_2^c\}\; >\; 0 \quad{\rm a.s.}$$
When $<A,\HH_2>\, =\r^n,$ one can then repeat a finite number of times the same discussion as above: alltogether this entails that $\tau > 0$ a.s. if $\nu^B$ does not exhaust $\r^n,$ which completes the proof by contraposition.

\qed

\subsection{Proof that $\nu^B$ exhausts $\r^n$ w.r.t. $A\,\Rightarrow\, X_1$ is a.c.} This is the difficult inclusion and we will first establish three lemmas. The first one is an easy application of the implicit function theorem. The second one is an a.s. independence result on a certain class of linear systems, for which we could not find any reference in the literature on control theory. The third one allows to choose suitably the jumping times of $B$ which will be later targeted into $X_1$ via a certain a.s. submersion. Throughout, $\lambda$ will stand for the Lebesgue measure independently of the underlying Euclidean space.

\begin{lemma}
\label{Inv1}
Let $X$ be an absolutely continuous random variable in $\r^p$ with $p\ge n$ and $\pph : \r^p \to \r^n$ be a $\CC^1$ function such that its Jacobian matrix $\d \pph$ verifies ${\rm Rank}\,\d \pph (X)  = n$ a.s. Then $\pph(X)$ is absolutely continuous in $\r^n.$  
\end{lemma}

\noindent
{\em Proof:} Choose any $x^{}_\NN\in\r^p\cap\NN^c,$ where $\NN = \{x\in\r^p\; /\; {\rm Rank}\;\d \pph (x)  < n\}.$ Setting 
$$\Xde\; =\; X\Un_{\{X\in\NN^c\}}\; +\; x^{}_\NN\Un_{\{X\in\NN\}},$$
we see by assumption that $X =\Xde$ a.s. so that $\pph(X) = \pph(\Xde)$ a.s. as well. Hence, it suffices to show that $\pph(\Xde)$ is absolutely continuous in $\r^n$. In the case $p = n$ this is a well-known fact for which we found a proof in \cite{St}, Lemma IV.3.1. For the sake of completeness, we will give an argument in the general case $p\ge n.$ 

Fix an underlying probability space $(\Omega, \FF,\p).$ By approximation, it is enough to show that for every relatively compact set $\Gamma\subset \Phi = \{\pph(\Xde(\omega)),\; \omega\in\Omega\}$ such that $\lla(\Gamma) =0,$ 
$$\p[ \pph(\Xde)\in\Gamma] \; =\; 0.$$ 
For every $y\in\Gamma,$ fix $x\in\pph^{-1}(y)\subset A.$ Since Rank$\d\pph(x)\, =\, n,$ by the implicit function theorem there exist $\VV_x$ and $\WW_y$ open neighbourhoods of $x$ and $y$ respectively, $\OO_p$ and $\OO_n$ open neighbourhoods of $0$ respectively in $\r^p$ and $\r^n$ endowed with a canonical basis, $\psi_x : \VV_x \to \OO_p$ and $\psi_y : \WW_y \to \OO_n$ diffeomorphisms such that 
$$\psi_y\,\circ\,\pph\,\circ\,\psi_x^{-1} \; : \; \OO_p \; \to \; \OO_n$$
is the canonical projection from $\OO_p$ to $\OO_n.$ Taking a finite covering $\lacc \WW_{y_1}, \ldots, \WW_{y_k}\racc$ of $\Gamma$ yields
\begin{eqnarray*}
\p[ \pph(\Xde)\in\Gamma] & \le & \sum_{i=1}^k \p [ \pph(\Xde) \in \WW_{y_i}\cap\Gamma]\\
& = & \sum_{i=1}^k \p [ \psi_{y_i}\circ \pph(\Xde) \in \Gamma_i]
\end{eqnarray*}   
with the notation $\Gamma_i = \psi_{y_i}(\WW_{y_i}\cap\Gamma),$ which has Lebesgue measure zero in $\r^n$ since $\psi_{y_i}$ is a diffeomorphism. Setting $\AA_i\; =\; \{\psi_{y_i}\circ \pph(\Xde) \in \Gamma_i\}$ for every $i\in\{1,\ldots,k\},$ the random variables
$Y_i = \Un_{\AA_i}\psi_{x_i}(\Xde)$ are absolutely continuous in $\r^p$ since $\psi_{x_i}$ are diffeomorphisms. Besides, since the projections $\psi_{y_i}\circ \pph \circ \psi_{x_i}$ have full rank, from Lemma 2.3 in \cite{PZ} the random variables $\psi_{y_i}\circ \pph \circ \psi_{x_i}(Y_i)$ are absolutely continuous in $\r^n$ . Hence, for every $i\in\{1,\ldots, k\},$
\begin{eqnarray*}
\p [\psi_{y_i}\circ \pph(\Xde) \in \Gamma_i] & = & \p \lcr \psi_{y_i}\circ \pph \circ \psi_{x_i}(Y_i) \in \Gamma_i\rcr\; = \; 0,
\end{eqnarray*}
which entails that $\p[ \pph(\Xde)\in\Gamma] \, =\, 0$ and completes the proof.

\qed

\begin{lemma} 
\label{Van1}
Assuming (\ref{Rank}), let $(b_1, \ldots, b_r)$ be a generating sequence with respect to $(A,B)$ for some $r\in[\k, m].$ Then the set
$$\lacc (t^1_1, \ldots, t_q^1, \ldots, t^r_1, \ldots, t_q^r)\in\r^{q\times r} \; /\;{\rm Rank} \lcr \ee^{t^1_1 A}b_1, \ldots, \ee^{t_q^1 A}b_1, \ldots, \ee^{t^r_1 A}b_r, \ldots, \ee^{t_q^r A}b_r\rcr\; <\; n\racc$$ 
has zero Lebesgue measure. 
\end{lemma}
\noindent
{\em Proof:} 
We first consider the case $\k = r = 1.$ Fix $b\in\BB$ a generating vector such that $\r^n = {\rm Vect} \lcr b, Ab,\ldots, A^{n-1}b\rcr.$ The function 
$$(t_1, \ldots, t_n)\;\mapsto\; {\rm Det} \lcr \ee^{t_1 A}b, \ldots, \ee^{t_n A}b\rcr$$
is analytic in $\r^n$ and it is not identically zero. Actually, if it were, then the analytic function $\rho : t\mapsto\ee^{tA}b$ would be valued in some fixed hyperplane of $\r^n,$ as well as all its successive derivatives, which is impossible since
$${\rm Vect} \lcr \rho(0), \rho'(0),\ldots, \rho^{(n-1)}(0)\rcr\; =\; {\rm Vect} \lcr b, Ab,\ldots, A^{n-1}b\rcr\; =\;\r^n.$$ 
Hence, since the zero set of a real analytic function over $\r^n$ either is $\r^n$ itself or has zero Lebesgue measure - see \cite{F} p. 240 or Lemma 2 in \cite{Y}, we obtain that
$${\rm Rank }\lcr \ee^{t_1 A}b, \ldots, \ee^{t_n A}b\rcr\; =\; n$$
almost everywhere in $\r^n$, which completes the proof when $\k=r=1.$

We now proceed to the remaining cases. Recalling (\ref{minime}), we first claim that the $q\times q$ matrix 
$$\Psi_q (t_1, \ldots, t_q)\; =\; \lacc \psi_i (t_j)\racc_{1\le i, j\le q}$$
has rank $q$ a.e. in $\r^q.$ Indeed, by the definition of $q$ there exists $b\in\r^n$ such that Rank $[b, Ab,\ldots, A^{q-1}b]\, =\, q.$ Setting $\AA_b\, =\, {\rm Vect} [b, Ab,\ldots, A^{q-1}b]$ and viewing $A$ as a cyclic endomorphism of the $q-$dimensional vector space $\AA_b,$ we see from the case $\k=r=1$ that Rank $[\ee^{t_1 A}b, \ldots, \ee^{t_q A}b]\; =\; q$ a.e. in $\r^q.$ However, it follows from (\ref{minime}) that
$$[\ee^{t_1 A}b, \ldots, \ee^{t_q A}b]\; =\; [b, \ldots, A^{q-1}b]\times\Psi_q (t_1, \ldots, t_q)$$
so that $\Psi_q (t_1, \ldots, t_q)$ must have rank $q$ a.e. in $\r^q$ as well. Let now $(b_1, \ldots, b_r)$ be a generating sequence with respect to $(A,B).$ Setting $\BB_i \, =\, {\rm Vect} [b_i, Ab_i,\ldots, A^{q-1}b_i],$ we have Dim $\BB_i\, \le \, q$ for every $i = 1\ldots r.$ Similarly as above, 
$$[\ee^{t_1 A}b_i, \ldots, \ee^{t_q A}b_i]\; =\; [b_i, \ldots, A^{q-1}b_i]\times\Psi_q(t_1, \ldots, t_q)$$ 
and since $\Psi_q$ is a.e. invertible, this entails that Rank $[\ee^{t_1 A}b_i, \ldots, \ee^{t_q A}b_i]\, =\,$ Dim $\BB_i$ a.e. in $\r^q.$ Notice that by (\ref{minime}), one has $\ee^{tA}b_i\in\BB_i$ for every $t\in\r,$ so that $[\ee^{t_1 A}b_i, \ldots, \ee^{t_q A}b_i]$ forms actually a basis of $\BB_i$ a.e. in $\r^q,$ for every $i= 1\ldots r.$ Besides, it follows from the definition of the Lebesgue measure that if $\AA_1, \ldots, \AA_r$ are negligible sets in $\r^q,$ then $(\AA_1^c\times\ldots\times\AA_r^c)^c$ is negligible in $\r^{q\times r}.$ Putting everything together entails that a.e. in $\r^{q\times r},$ 
$${\rm Rank} \lcr \ee^{t^1_1 A}b_1, \ldots, \ee^{t_q^1 A}b_1, \ldots, \ee^{t^r_1 A}b_r, \ldots, \ee^{t_q^r A}b_\k\rcr\; =\; {\rm Dim}(\BB_1 + \ldots + \BB_r)\; =\; n,$$
where the last equality comes from the definition of the generating sequence $(b_1, \ldots, b_r).$ The proof is complete.

\qed

\begin{remarks} {\em (a) By the same argument, one can prove that 
$$\lla \lacc (t_1, \ldots, t_n)\in\r^{n} \; /\;{\rm Rank} \lcr \ee^{t_1 A}B, \ldots, \ee^{t_n A}B\rcr\; <\; n\racc\; =\; 0$$
as soon as (\ref{Rank}) holds. An interesting point is that when $A$ has real spectrum, then (\ref{Rank}) actually entails that 
\begin{equation}
\label{exact}
{\rm Rank} \lcr \ee^{t_1 A}B, \ldots, \ee^{t_n A}B\rcr\; =\; n
\end{equation}
for {\em every} $\,t_1 < \ldots < t_n.$ The latter is false nevertheless when $A$ has non-real eigenvalues. The proofs of the above two facts involve rather technical considerations on generalized Vandermonde matrices, which we shall not discuss here. 

\vspace{1mm}

\noindent
(b) I could not find in the literature any material on the following question, whose positive answer would quickly entail (\ref{exact}). Assuming that $A$ has real spectrum and that (\ref{Rank}) holds, let $u^1, \ldots, u^n$ be $n$ piecewise constant control functions which are linearly independent over $[0,1]$ and consider the linear systems
$$\frac{d x^i_t}{d t}\; = \; x\; +\; \int_0^t (Ax^i_s + Bu^i_s)\, \d s, \quad t\ge 0, \;\; i =1\ldots n.$$
Does $(x^1_1, \ldots, x^n_1)$ form then necessarily a basis of $\r^n$?
}
\end{remarks}

A family $(\CC_1, \ldots, \CC_r)$ of disjoint pointed cones with common vertex at zero such that every $(c_1, \ldots, c_n)\in(\CC_1, \ldots, \CC_r)$ is a generating sequence with respect to $(A,B)$ will be called a {\em generating garland} with respect to $(A,B),$ or simply a generating garland when there is no ambiguity. When $(b_1, \ldots, b_r)$ is a generating sequence, notice that $(\CC_1, \ldots, \CC_r)$ defined by $\CC_i =\{\mu b_i, \; \mu >0\}$ for $i=1\ldots r,$ is a generating garland. Our last lemma makes a connection between this notion and the exhaustion property:
 
\begin{lemma}
\label{Geom}
If $\nu^B$ exhausts $\r^n$ w.r.t. $A$, then for every $M>0$ there exists $r\in[\kappa,m]$ and a generating garland $(\CC_1, \ldots, \CC_r)$ such that $\nu^B(\CC_i) \ge M$ for every $i = 1\ldots r.$
\end{lemma}

\noindent
{\em Proof.} Fix $M > 0.$ If $\nu^B$ exhausts $\r^n$ w.r.t. $A$, then we know that $<A,\BB> \, =\r^n.$ Suppose first that $\nu(\HH^c) =\infty$ for every hyperplane $\HH\subset \BB$ and let us show that there exists a generating garland $(\CC_1, \ldots, \CC_m)$ such that $\nu^B(\CC_i) \ge M$ for every $i = 1\ldots m,$ which is intuitively obvious. 

If $\SS^{m-1}$ denotes the unit Euclidean sphere of $\BB,$ consider a family $\{\Pi_\delta, \; \delta > 0\}$ of finite measurable partitions of $\SS^{m-1}$ such that ${\rm Diam} \,(\Pi_\delta)\, < \, \delta$ and $\Pi_{\delta'}$ is a subpartition of $\Pi_\delta$ for every $\delta' < \delta.$ Let $\CC^\delta_M$ denote the disjoint finite family of pointed cones with vertex at zero and apex in $\Pi_\delta$ such that $\nu^B(\CC) \ge M$ for every $\CC\in\CC^\delta_M. $ If for every $\delta > 0$ no generating garland of size $m$ is contained in $\CC^\delta_M,$ then for every $\delta > 0$ there exists at least one hyperplane $\HH_\delta$ intersecting every $\CC\in\CC^\delta_M,$ and the assumption ${\rm Diam} \,(\Pi_\delta)\, < \, \delta$ readily entails that all these $\HH_\delta$'s converge - in the sense that their normal unit vectors converge in the metric space $\SS^{m-1}$ - to some fixed hyperplane $\HH_0.$ Last, it is a bit tedious but not difficult to see that by construction and by the finiteness of $\Pi_{\delta},$ one must have $\nu(\HH_0^c) <\infty,$ which yields a contradiction. Hence, there exists $\delta > 0$ such that $\CC^\delta_M$ contains a generating garland of size $m,$ and we are done.

If $\nu(\HH^c) <\infty$ for some hyperplane $\HH\subset \BB$, then by definition of the exhaustion property there must exist an hyperplane $\HH_{m-1}$ such that $<A,\HH_{m-1}> \, =\r^n.$ If $\nu(\HH_{m-1}\cap\HH^c) =\infty$ for every subspace $\HH\subset \HH_{m-1},$ then reasoning exactly as above we can show that there exists a generating garland $(\CC_1, \ldots, \CC_{m-1})$ such that $\nu^B(\CC_i) \ge M$ for every $i = 1\ldots m-1.$ If $\nu(\HH_{m-1}\cap\HH^c) <\infty$ for some subspace $\HH\subset \HH_{m-1},$ then we can repeat the same procedure as above. But again by the definition of the exhaustion property, the latter procedure cannot be repeated more than $(m-\kappa)$ times: all in all, this shows that there exists $r\in[\kappa,m]$ and a generating garland $(\CC_1, \ldots, \CC_r)$ such that $\nu^B(\CC_i) \ge M$ for every $i = 1\ldots r.$

\qed

\begin{remark}{\em When $\nu^B$ exhausts $\r^n,$ there may exist no generating garland $(\CC_1, \ldots, \CC_r)$ such that $\nu^B(\CC_i) = \infty $ for every $i = 1\ldots r.$ Think of the situation where $\kappa = m = 2$ and $\nu^B$ is infinite with support in the arc $y = x^2$, $\BB$ being endowed with an orthonomal frame $Oxy$.}
\end{remark}

\noindent
{\bf End of the proof.} From (\ref{levy}) and after time-reversal, it is enough to prove that
$$Y\; =\; \int_0^1\ee^{sA}B \d Z_s\; =\; \int_0^1\ee^{sA}\d B_s$$
is absolutely continuous. For this purpose, we will use the same method as depicted in the introduction, in a somewhat more elaborated manner. If $\Gamma\subset \r^n$ is such that $\lla(\Gamma) = 0,$ we need to show that for every $\eps > 0$ 
$$\p[Y\in\Gamma]\; < \; \eps.$$
Fix $\eps >0$ and let $M >0$ be such that
$$\p[T^M_{q+1}\ge 1]\; <\; \eps/m$$
where $T^M_{q+1}$ is the sum of $(q+1)$ independent exponential variables with parameter $M$. By Lemma \ref{Geom}, there exists $r\in[\kappa,m]$ and a generating garland $(\CC_1, \ldots, \CC_r)$ such that $\nu^B(\CC_i) \ge 2M$ for every $i = 1\ldots r.$ Besides, if $\BB_\eta$ stands for the Euclidean ball of $\BB$  centered at 0 with radius $\eta$ and if $\CC_i^\eta = \CC_i\cap \BB_\eta^c,$ then we can actually choose $\eta > 0$ such that $\nu^B(\CC_i^\eta) \ge M$ for every $i = 1 \ldots r.$ Let $\{ T_{i,p}^\eta, \; p\ge 1\}$ be the ordered sequence of jumping times of the L\'evy process $\{B_t, \; t\ge 0\}$ into $\CC^\eta_i$ and set $T_p^\eta = \sup\{T_{i,p}^\eta, \; i =1\ldots r\}$ for every $p\ge 1.$ We have
$$\p[T^\eta_{q+1}\ge 1]\; \le\; \sum_{i=1}^r \p[T^\eta_{i, q+1}\ge 1]\;\le\;r \p[T^M_{q+1}\ge 1]\; <\; r\eps/m\; \le \; \eps$$
and it is hence sufficient to prove that
\begin{equation}
\label{tite}
\p\lcr T_{q+1}^\eta < 1, \; Y\in\Gamma\rcr\; =\; 0
\end{equation}
for every $\eta > 0.$ Let $\FF^\eta$ be the $\sigma-$algebra generated by $\lacc T_{i,p}^\eta, \; p\geq q+1, \; i = 1\ldots r\racc,$ $\{ \Delta B_{T_{i,p}^\eta}, \; p\geq 1, \; i = 1\ldots r\}$ and the L\'evy process ${\tilde B}^\eta$ defined by
$${\tilde B}^\eta_t\; = \; B_t - \sum_{T_{i,p}^\eta\leq t} \Delta B_{T_{i,p}^\eta}, \quad t\geq 0.$$
On $\{ T_{q+1}^\eta < 1\},$ one can write
$$Y\; =\; Y^\eta\; +\; \sum_{i=1}^{r}\sum_{j=1}^q\ee^{T_{i,j}^\eta A}\Delta B_{T_{i,j}^\eta}$$
with $Y^\eta$ a $\FF^\eta-$measurable random variable. Since $T_{q+1}^\eta$ is $\FF^\eta-$measurable as well, we have
\begin{eqnarray*}
\p\lcr T_{q+1}^\eta < 1, \; Y\in\Gamma\rcr & = & \p\lcr T_{q+1}^\eta < 1, \; \p\lcr Y\in\Gamma\;\vert\;\FF^\eta\rcr\rcr\\
& = & \p[T_{q+1}^\eta < 1, \; \p[ \Ydee\in\Gamma^\eta\;\vert\;\FF^\eta]]
\end{eqnarray*}
where $\Gamma^\eta = \Gamma - Y^\eta$ is a $\FF^\eta-$measurable set such that $\lla(\Gamma^\eta) = 0,$ and
$$\Ydee\; =\;\sum_{i=1}^{r}\sum_{j=1}^q\ee^{T_{i,j}^\eta A}\Delta B_{T_{i,j}^\eta}.$$
The key-point is that by standard properties of jumping measures, since the $\CC_i^\eta$'s are disjoint, conditionally on $\FF^\eta$ the law of the $(r\times q)-$uple $(T_{1,1}^\eta,\ldots ,T_{1,q}^\eta, \ldots, T_{r,1}^\eta,\ldots ,T_{r,q}^\eta)$ is that of the tensor product of $r$ independent $q$-th order statistics respectively on $[0, T_{i, q+1}^\eta],$ viz. the tensor product of $r$ independent uniform laws on the respective sets
$$\lacc 0 < t^i_1 < \ldots < t^i_q < T_{i, q+1}^\eta\racc.$$
In particular, the law of $(T_{1,1}^\eta,\ldots ,T_{1,q}^\eta, \ldots, T_{r,1}^\eta,\ldots ,T_{r,q}^\eta)$ is absolutely continuous in $\r^{q\times r}$ and by Lemma \ref{Inv1}, (\ref{tite}) will hold as soon as the Jacobian matrix of the application 
$$(t^1_1,\ldots ,t^1_q, \ldots, t^r_1,\ldots ,t^r_q)\; \mapsto\; \sum_{i=1}^{r}\sum_{j=1}^q\ee^{t^i_j A}\Delta B_{T_{i,j}^\eta}$$
from $\r^{q\times r}$ to $\r^n$ has rank $n$ a.e. conditionnally on $\FF^\eta$ (recall that the sequence $\{ \Delta B_{T_{i,p}^\eta}, \; p\geq 1, \; i = 1\ldots r\}$ is $\FF^\eta$-measurable and independent of $\{ T_{i,p}^\eta, \; p\leq q, \; i = 1\ldots r\}$). This Jacobian matrix is equal to
\begin{equation}
\label{yakov}
A\times\lcr \ee^{t^1_1A}\Delta B_{T_{1,1}^\eta},\ldots, \ee^{t^1_q A}\Delta B_{T_{1,q}^\eta}, \ldots, \ee^{t^r_1A}\Delta B_{T_{r,1}^\eta},\ldots, \ee^{t^r_q A}\Delta B_{T_{r,q}^\eta}\rcr
\end{equation}
and, since $A$ is invertible and $\Delta Z^{B}_{T_{i,j}^\eta}\in \CC_i^\eta$ for every $i=1\ldots r$ and $j=1\ldots q,$ conditionnally on $\FF^\eta$ it has a.s. the same rank as
$$\lcr \ee^{t^1_1 A}b_1, \ldots, \ee^{t_q^1 A}b_1, \ldots, \ee^{t^r_1 A}b_r, \ldots, \ee^{t_q^r A}b_r\rcr$$
where $(b_1, \ldots, b_r)$ is some generating sequence of $\r^n$ with respect to $(A,B).$ Now by Lemma \ref{Van1}, the latter has full rank a.e. and the proof is finished.

\qed

\begin{remark}{\em The invertibility assumption on $A$ is only useful for getting the full rank a.s. of the Jacobian matrix given in (\ref{yakov}). Nevertheless this is a crucial assumption, and in the next section we will give a counterexample when $A$ is singular.}
\end{remark}
\section{Final remarks}

\subsection{The case with a Brownian component} If the driving L\'evy process
$Z$ has a non-trivial Gaussian component, then by the linearity of (\ref{OU2})
this amounts to consider the problem of absolute continuity for
$$X_t \; = \; \ee^{tA}x\; +\; \int_0^t \ee^{(t-s)A}\d W_s\; +\; \int_0^t \ee^{(t-s)A}\d B_s, \quad t\ge 0,$$
where $W$ is some $\BB$-valued Brownian motion independent of the Non-Gaussian
L\'evy process $B$. Set $H = \,<A, {\rm Im}\, W>$ and, given any Euclidean
structure on $\r^n,$ denote by $H^\perp$ the orthogonal complement of $H$ in
$\r^n.$ The $H-$valued random variable
$$\int_0^1 \ee^{(1-s)A}\d W_s$$ 
is Gaussian and by a classical result in control theory - see e.g. Theorem 1.1
in \cite{W} - it is non-degenerated, hence absolutely continuous in
$H$. Since $W$ and $B$ are independent, Lemma 3 in \cite{Y} and Lemma 2.3
in \cite{PZ} yield
$$X_1\; {\rm is \; a.c.}\;\; \Longleftrightarrow\;\; \Pi_{H^\perp}\lpa\int_0^1
\ee^{(1-s)A}\d B_s\rpa\; {\rm is \; a.c.}$$
where $\Pi_{H^\perp}$ stands for the orthogonal projection operator onto
$H^\perp.$ A straightforward modification of our proof entails then
$$X_1\; {\rm is \; a.c.}\;\; \Longleftrightarrow\;\; \tau^{}_H\, =\, 0\;\, {\rm
  a.s.}\;\; \Longleftrightarrow\;\;  \nu^B \;{\rm exhausts}\; \r^n\; {\rm
  w.r.t.}\; (A,H).$$ 
where with the notations of the introduction we set $\tau^{}_H\; =\; \inf\{ t>0, \; \AA_t +
H \, =\r^n\},\k^{}_H$ for the minimal number of linearly independent vectors
$b_1, \ldots,b_p\in\BB$ such that  
$$<A,b_1>\, +\ldots + <A,b_p> \,+ \; H  =\r^n,$$ 
and say that $\nu^B$ exhausts $\r^n$ with respect to $(A,H)$ if $<A,\BB> \,+ \; H  =\r^n$ and
there exists $r\in [\k^{}_H, m]$ and a subspace $\HH_r\subset\BB$ of
dimension $r$ such that $<A,\HH_r>\, + \; H  =\r^n$ and $\nu^B(\HH_r\cap\HH^c)
=+\infty$ for every hyperplane $\HH\subset\HH_r.$

\subsection{Some explicit descriptions of the exhaustion property.}
From now on we will assume that $Z$ has no Gaussian part, that is $B$ has no
Gaussian part either. Let us first consider the case $n=1,$ i.e. $X$ is solution to 
\begin{equation}
\label{OU1}
\d X_t\; =\; aX_t  \d t \; +\; \d B_t
\end{equation}
where $a\in\r$ and $B$ is one-dimensional. The exhaustion property just means that $\nu^B$ is infinite and our result reads
\begin{equation}
\label{n1}
X_1\; {\rm is \; a.c.}\;\; \Longleftrightarrow\;\; \nu^B \;{\rm is\; infinite}
\end{equation}
as soon as $a\neq 0.$ Notice that this is actually an immediate consequence of
Theorem A in \cite{NS} - see also Theorem 1.1. in \cite{K}. Let us also give a short proof of the non-trivial reverse inclusion in (\ref{n1}), similar to that given in the introduction: the solution to (\ref{OU1}) is given by
$$X_t \; = \; \ee^{ta}x\; +\; \int_0^t \ee^{a(t-s)}\d B_s\; = \; \ee^{ta}x\; +\; B_t \; + \; a\int_0^t \ee^{a(t-s)} B_s \d s, \quad t\ge 0,$$
where $x$ is the initial condition and where in the second equality we made an
integration by parts, assuming $B_0 =0$ without lost of generality. Hence, leaving the details to the reader, one may follow roughly the same method as for the integral of $B$, noticing with the same notations that on $\{T^\eta_2 < 1\}$ the value of $B_1$ does not depend on $T_1^\eta$, hence $B_1$ is $\FF^\eta$-measurable as well. \\

Let us now discuss the case $n=2.$ To simplify the notations, we will denote
by $\II$ the set or family of sets where $\nu^B$ is infinite if and only if
the exhaustion property holds. We will also suppose implicitly that $A$ is non-singular. Up to some equivalent transformations on $A$ which are not relevant to the absolute continuity problem, there are four situations:

\vspace{1mm}

\noindent
(a) $A$ has no real eigenvalue, in other words $A$ is a multiple of
$$\lpa\begin{array}{rl} \cos\theta & \sin\theta \\
-\sin\theta & \cos\theta
\end{array}\rpa$$
for some $\theta\in]0,\pi[.$ Then $\k = 1$ i.e. $A$ is cyclic, and it is easy to see that every non-zero vector of $\r^2$ is generating. Hence we simply have $\II = \r^2$ viz.  as in the real case, $X_1$ is a.c. if and only if $\nu^B$ is infinite.

\vspace{1mm}

\noindent
(b) $A$ is a multiple of the identity matrix. Then $\k = 2$ and $<A,b> \, =
{\rm Vect}\{b\}$ for every $b\in\r^2.$ This means that $\II =\{({\rm
  Vect}\{b\})^c, \; b\in\r^2\}$ viz. the infinite part of $\nu^B$ must not be carried by any line in $\r^2$.

\vspace{1mm}

\noindent
(c) $A$ is a Jordan cell matrix, i.e. $A$ is of the type
$$\lpa\begin{array}{cc} \alpha & 1 \\
0 & \alpha
\end{array}\rpa$$
with $\alpha\neq 0.$ Then $\k = 1$ and every non-zero vector of $\r^2$ is generating except the multiples of $(1,0):$ we have $\II = ({\rm Vect}\{(1,0\})^c.$

\vspace{1mm}

\noindent
(d) $A = \mbox{Diag}(\alpha, \beta)$ with $\alpha\neq\beta$ and $\alpha, \beta$ non zero. Then $\k = 1$ and every non-zero vector of $\r^2$ is generating except those in Vect $\{(1,0)\}\cup\mbox{Vect}\{(0,1)\}.$ On the other hand, (Vect $\{(1,0)\} - \{0\}$, Vect $\{(0,1)\}-\{0\})$ is a generating garland: we have 
$\II = \{({\rm Vect}\{(1,0)\})^c\cap({\rm Vect}\{(0,1)\})^c\} \cup \{{\rm Vect}\{(1,0)\} \;\mbox{and Vect}\{(0,1)\}\}.$ \\

When $n > 2,$ it becomes quite lengthy to depict the exhaustion property. Let us give however four typical examples when $n=3,$ keeping for $\II$ the same meaning as above and using the notations $\HH_x = \{x=0\}, \HH_y =  \{ y=0\}$ and $\HH_z = \{z=0\}$ where $Oxyz$ is a given orthogonal frame of $\r^3.$ 

\vspace{1mm}

\noindent
(f) $A$ is a Jordan cell matrix, i.e. $A$ is of the type
$$\lpa\begin{array}{ccc} \alpha & 1 & 0\\
0 & \alpha & 1\\
0 & 0 & \alpha
\end{array}\rpa$$
with $\alpha\neq 0.$ Then $\k = 1$ and every non-zero vector of $\r^3$ is generating except those in $\HH_z:$  we have $\II = \HH_z^c.$

\vspace{1mm}

\noindent
(g) $A$ is a block matrix of the following type
$$\lpa\begin{array}{ccc} \alpha & 0 & 0\\
0 & \beta & 1\\
0 & 0 & \beta
\end{array}\rpa$$
with $\alpha\neq\beta$ and $\alpha, \beta$ non zero. Then $\k = 1$ and every non-zero vector of $\r^3$ is generating except those in $\HH_x\cup\HH_z:$  we have $\II = \HH_x^c\cap\HH_z^c.$

\vspace{1mm}

\noindent
(h) $A$ is a block matrix of the following type
$$\lpa\begin{array}{ccc} \alpha & 0 & 0\\
0 & \alpha & 1\\
0 & 0 & \alpha
\end{array}\rpa$$
with $\alpha\neq 0.$ Then $\k = 2$ and every generating sequence must not be
valued in any hyperplane $\HH_u = u^\perp$ with $u$ a unit vector of $Oxz:$ we have $\II = \{\HH_u^c, \, u\in Oxz\}.$

\vspace{1mm}

\noindent
(i) $A = \mbox{Diag}(\alpha, \beta, \gamma)$ with distinct non zero $\alpha,
\beta$ and $\gamma.$ Then $\k = 1$ and every vector in
$\HH^c_x\cap\HH_y^c\cap\HH_z^c$ is generating. But as in dimension 2, one can
also build generating garlands with one component in $\HH_x, \HH_y$ or
$\HH_z.$ The infinity set $\II$ is then the union of the following eight sets
or families of sets: $\HH_x^c\cap\HH_y^c\cap\HH_z^c, \{\HH_x\cap\HH_y^c\cap\HH_z^c\,\mbox{and}\,\HH_y\cap\HH_x^c\}, \{\HH_y\cap\HH_x^c\cap\HH_z^c\,\mbox{and}\,\HH_x\cap\HH_y^c\}, \{\HH_y\cap\HH_z^c\cap\HH_x^c\,\mbox{and}\,\HH_z\cap\HH_y^c\}, \{\HH_z\cap\HH_y^c\cap\HH_z^c\,\mbox{and}\,\HH_y\cap\HH_z^c\}, \{\HH_z\cap\HH_x^c\cap\HH_y^c\,\mbox{and}\,\HH_x\cap\HH_z^c\}, \{\HH_x\cap\HH_z^c\cap\HH_y^c\,\mbox{and}\,\HH_z\cap\HH_x^c\}, \{\HH_x\cap\HH_y\cap\HH_z^c\,\mbox{and}\,\HH_y\cap\HH_z\cap\HH_x^c\,\mbox{and}\,\HH_z\cap\HH_x\cap\HH_y^c\}.$ 

\subsection{Some open questions} As we mentioned before, our theorem no more
holds when $A$ is singular, as shows the following counterexample with $n=2,d=m=1,$
$$A\; =\, \lpa\begin{array}{ll} 0 & 0 \\
1 & 0 \end{array}\rpa\quad\mbox{and}\quad B\; =\, \lpa\begin{array}{l} 1 \\
0 \end{array}\rpa.$$
From the control theory viewpoint, this example yields the so-called rocket car equations, which
serve as toy-models  \cite{JM} for studying the Pontryagin maximum
principle. From the stochastic
viewpoint, assuming $x =(0,0)$ one gets the so-called Kolmogorov process
$X = (X^1, X^2),$ with
$$X^1_t\; =\; Z_t\quad\mbox{and}\quad X^2_t\; =\; \int_0^t Z_s \d s,\quad t\ge
0,$$
and $Z$ is a one-dimensional L\'evy process. Notice that $(A,B)$ is
controllable, so that $\nu^B$ exhausts $\r^2$ w.r.t. $A$ if and only if $\nu^Z$ is
infinite. But then $X^1_1 = Z_1$ may be singular - see again Theorem 27.19 in
\cite{S} - so that by Lemma 2.3 in \cite{PZ}, $X_1$ is not absolutely
continuous either. 

When $A$ is singular and $m$, one may wonder if the following holds true for
every $t > 0:$
$$B_t\; {\rm is \; a.c. \; in}\;\r^n\; \Longrightarrow\;\; X_t\; {\rm is \;
  a.c.\; in}\; \r^n.$$
From our theorem, this property is trivial when $n = 1$ and we also refer to
Theorem B in \cite{NS} for a non-linear extension. When $n\ge 2$ some
geometrical difficulties arise however, which will be the matter of further
research. When $m <n$ the problem seems much more complicated without further
conditions on the jumping measure. In particular I do not have the answer to
the following basic question, which would solve the problem at least for
the Kolmogorov process:
\begin{center}
For a real L\'evy process $\{Z_t, \, t\ge 0\},$ if $Z_1$ is a.c. in $\r,$ is ${\displaystyle \lpa Z_1, \int_0^1Z_t \d t\rpa}$
a.c. in $\r^2$?
\end{center} 

To conclude this paper, let us go back to the case $n=1$ and consider the
following class of infinitely distributions
$$\OO\; =\; \lacc\LL\lpa \int_0^1 \ee^s \d B_s\rpa, \; B\;\mbox{real L\'evy process}\racc,$$
for which we proposed the name OU class. If $\mu^B\in\OO$ corresponds to some
Non-Gaussian L\'evy process $B$, after time-reversal our previous discussion
entails 
\begin{equation}
\label{class}
\mu^B\; {\rm is \; a.c.}\;\; \Longleftrightarrow\;\; \nu^B \;{\rm is\; infinite.}
\end{equation}
Besides, with a little reflexion, one can show that (\ref{class}) also holds when replacing in $\OO$ the kernel $\ee^s$ by any $\CC^1$ function $f(s)$ whose derivative does not vanish in $]0,1[.$ In particular, the equivalence (\ref{class}) will hold for the well-known U-class where $f(s) =s,$ and B-class where $f(s) = \log s.$ Of course, for these two special classes one could get (\ref{class}) directly in considering the jumping measure of $\mu^B$ which happens to be absolutely continuous - see \cite{BMS} for details - and applying Tucker's result - see Theorem 27.7 in \cite{S}. It would be interesting to investigate which exact class of integration kernels entails (\ref{class}) for L\'evy integrals, and also what occurs in the multivariate case.

\end{document}